\documentclass[12pt]{amsart}

\usepackage{latexsym}
\usepackage{amsmath}
\usepackage{amssymb}
\usepackage{a4}

\newcommand{\B}{{\mathbf B}}

\newcommand{\R}{{\mathbf R}}

\newtheorem{thm}{Theorem}[section] 
\newtheorem{prop}[thm]{Proposition} 
\newtheorem{lem}[thm]{Lemma} 
 
\newtheorem{cor}[thm]{Corollary} 
\newtheoremstyle{mydef}{}{}{\rmfamily}{}{\rmfamily}{.}{ }{\noindent\textbf{\thmname{#1}\thmnumber{ #2}\thmnote{ (#3)}}}
\theoremstyle{mydef}
\newtheorem{exa}[thm]{Example} 
\newtheorem{de}[thm]{Definition}
\numberwithin{equation}{section}

\newcommand{\dev}{\mathrm{dev}}

\newcommand{\ord}{\mathrm{ord}}

\newcommand{\GF}{\mathrm{GF}}
\newcommand{\SL}{\mathrm{SL}}

\newcommand{\0}{{\mathbf 0}}
\newcommand{\1}{{\mathbf 1}}
\newcommand{\2}{{\mathbf 2}}

\newcommand{\fixqed}{\qed\vskip+0.3\baselineskip}

\title
{Difference methods and Ferrero pairs}
\author{Tim Boykett and Peter Mayr}
\date{\today}
\thanks{AMS classification 16Y30. \\ \indent
 This work has been supported by grant P15691 of the Austrian National Science Foundation (Fonds zur F\"orderung der wissenschaftlichen Forschung).}

\begin{document}
\maketitle
\begin{center}
\begin{small}
Institut f\"ur Algebra, Stochastik und \\
wissensbasierte mathematische Systeme, \\
Johannes Kepler Universit\"at Linz, \\
4040 Linz,
Austria \\
{\tt tim@algebra.uni-linz.ac.at}\\
{\tt peter.mayr@algebra.uni-linz.ac.at}\\[5mm]
\end{small}
\end{center}

\begin{abstract}
 We present a 
 construction method of BIB-designs from a finite group $G$ and a group of
 automorphisms $\Phi$ on $G$ such that $|\Phi(x)| = |\Phi|$ for all $x\in G$,
 $x\neq 0$.
 By using a generalization of the concept of a difference family we can so
 unify several previous constructions of BIB-designs from planar near-rings.

\end{abstract}

\section{Introduction}

  This paper introduces a general construction that encompasses 
several known constructions of designs from near-rings.
The original work of Ferrero, lying historically between that
of Bose and Wilson, introduced the construction, related to
difference families. This paper uses a slightly more
general concept, the short difference family, to
generalize these constructions.

We proceed first by introducing our terminology, then the
short difference family construct.
We show that such constructions arise naturally from fixed-point-free
 automorphism groups acting on groups.

\section{Designs and Difference Families}

 Let $V$ be a set of size $v$ with $v>1$. Let $k,\lambda$ be positive integers.
 For our purposes, a $(v,k,\lambda)$-\emph{BIB-design} (or \emph{design} for 
 short) is a set $\B$ of subsets of $V$  that satisfies the following
 conditions: 
\begin{enumerate}
\item $|B|=k$ for all $B\in\B$; 
\item $|\{ B\in\B\ |\ u,v\in B\}|=\lambda$ for all $u,v\in V, u\neq v$ . 
\end{enumerate}
 The elements of $V$ are usually referred to as points; the elements of $\B$
 are called blocks.
 By our definition a design does not have repeated (multiple) blocks.
 Apart from that our definition is equivalent to that of a $2$-balanced design
 \cite[p. 15]{Beth:DesignI} or $(v,k,\lambda)$-BIBD \cite{Colbourn:Designs}.

 An \emph{automorphism} of a design $\B$ is a bijection of $V$ that fixes
 $\B$. 



 In this note we will use the following concepts:
 Let $(G,+)$ be a finite group, not necessarily abelian, of order $v>1$.
 We denote $G\setminus\{0\}$ by $G^*$.
 The group $G$ acts on itself by right translation.
 For $B\subseteq G$, we let $G_B = \{ g\in G\ |\ B+g = B\}$ denote the
 stabilizer of the set $B$ under the action of $G$.

 Let $\B$ be a set of subsets of $G$.
 Then $\dev\B = \{ B+g\ |\ B\in\B, g\in G\}$ is called the \emph{development}
 of $\B$ in $G$.
 For $B,C\in\B$, we define
\[ B\sim C \text{ if and only if } \exists g\in G: B = C+g. \]
 This relation $\sim$ is an equivalence relation on $\B$.
 The equivalence class of $B\in\B$ modulo $\sim$ is denoted by $B/\sim$.

\begin{de} \label{def:sdf}
 Let $\B$ be a set of subsets of the group $(G,+)$ with $|G|>1$.
 Let $v=|G|$, and let $k,\mu,\nu,\lambda'$ be positive integers such that the
 following are satisfied:
\begin{enumerate}
\item \label{it:sdf1}
 $|B| = k$ for all $B\in\B$;
\item \label{it:sdf2}
 $|G_B| = \mu$ for all $B\in \B$;
\item \label{it:sdf3} 
 $|B/\sim| = \nu$  for all $B\in \B$;
\item \label{it:sdf4}
 $|\{(B,a,b)\ |\ B\in\B, a,b\in B, a-b=d\}|=\lambda'$ for all 
 $d\in G, d\neq 0$. 
\end{enumerate}
 Then $\B$ is called a
 $(v,k,\frac{\lambda'}{\mu\nu})$-\emph{short difference family}
 (or \emph{sdf}) in $G$.
\end{de} 

 We note that an sdf with $\mu = 1$ is in fact a $(v,k,\lambda)$-difference
 family as defined in \cite[p. 470]{Beth:DesignI}.
 However, the general definition of a difference family allows multiple blocks,
 which we avoid.
 The term short difference family seems appropriate since in the literature
 a block $B\in\B$ with $|G_B|>1$ is said to have a short orbit.

 In \cite[p.271]{Colbourn:Designs} a set of subsets $\B$ of $G$ is called
 a partial difference family if $\dev\B$ is a BIB-design. The following
 proposition says that short difference families are partial difference
 families in this sense.

\begin{prop} \label{sdf}
 Let $\B$ be a $(v,k,\lambda)$-sdf in the group $G$. Then $\dev\B$ is a
 $(v,k,\lambda)$-design.
\end{prop}

 We note that $G$ acts as group of design-automorphisms on $\dev\B$ for an
 sdf $\B$. Hence the automorphism group of $\dev\B$ acts transitively on 
 the point set $G$ of $\dev\B$.

\emph{Proof of Proposition~\ref{sdf}.}
 Let $\B$ be a $(v,k,\lambda)$-sdf. 
 Then there exist $\mu,\nu$ such that $|G_B| = \mu$ and $|B/\sim| = \nu$
 for all $B\in \B$.
 We have $|\{ (B,a,b)\ |\ B\in\B,a,b\in B, a-b = d \}| = \lambda\mu\nu$
 for all $d\in G, d\neq 0$.

 If $B\sim C$ for $B,C\in\B$, then $\dev\{B\} = \dev\{C\}$. 
 Let $\R$ be a set of representatives of $\sim$ on $\B$. Then
 \[ \dev\B = \bigcup_{B\in\R} \dev\{B\}, \]
 where the union is disjoint. 
 
 Let $u,v\in G, u\neq v$. We count the number of blocks in $\dev\B$ that
 contain $\{u,v\}$:
\begin{eqnarray*}
 |\{C\in\dev\B\ |\ u,v\in C\}| & = & \sum_{B\in\R} |\{C\in\dev\{B\}\ |\ u,v\in C\}| \\
     & = & \sum_{B\in\B} \frac{|\{C\in\dev\{B\}\ |\ u,v\in C\}|}{|B/\sim|} \\
     & = & \frac{1}{\nu}\sum_{B\in\B} \frac{|\{ g\in G\ |\ u,v\in B+g\}|}{|G_B|} \\
     & = & \frac{1}{\mu\nu} \sum_{B\in\B} |\{ g\in G\ |\ u-g,v-g\in B\}|
\end{eqnarray*}
 For $B\in\B$ fixed, we now consider the map
 $f:\{ g\in G\ |\ u-g,v-g\in B\}\rightarrow\{(a,b)\in B\times B\ |\ a-b=u-v\}$,
 $g\mapsto(u-g,v-g)$. Then $f$ has an inverse, given by
 $f^{-1}:(a,b)\mapsto -a+u$. In particular, $f$ is bijective. Hence we obtain
\begin{eqnarray*}
 |\{C\in\dev\B\ |\ u,v\in C\}| & = & \frac{1}{\mu\nu} \sum_{B\in\B} |\{(a,b)\in B\times B\ |\ a-b=u-v\}| \\
     & = & \frac{1}{\mu\nu} \cdot |\{ (B,a,b)\ |\ B\in\B,a,b\in B, a-b = u-v \}| \\
     & = & \frac{1}{\mu\nu} \cdot \lambda\mu\nu \\
     & = & \lambda.
\end{eqnarray*}
 Thus each set $\{u,v\}$ with $u,v\in G, u\neq v$, is contained in
 $\lambda$ distinct blocks of $\dev\B$.  
 This completes the proof that $\dev\B$ is a $(v,k,\lambda)$-design.
\fixqed

 In the proof above we used that a set of representatives $\R$ for $\sim$ on
 $\B$ has the same development as $\B$.
 Now each equivalence class of $\sim$ on $R$ has size $1$.
 Thus we could have required $|B/\sim| = 1$ for all $B\in\B$ instead of
 condition~\eqref{it:sdf3} in our definition of an sdf, and we would still
 obtain the same designs.
 It is just a matter of notational convenience that we allow $|B/\sim| = \nu$
 for an arbitrary, but fixed integer $\nu$.

 In the following lemma $G_B$ is characterized as the unique maximal subgroup
 of which $B\subseteq G$ can be expressed as union of left cosets.

\begin{lem} \label{le:G_B}
 Let $G$ be a group, and let $B\subseteq G$. Then $G_B$ is the unique maximal
 (with respect to containment) element in $\{H\leq G\ |\ B = B+H\}$.
\end{lem}

\begin{proof} Straightforward.
%
\end{proof}

 As a corollary of Lemma~\ref{le:G_B}, we obtain that $G_B=\{0\}$ if
 $|B|$ and $|G|$ are relatively prime.

\section{Short difference families from fpf automorphisms}

 Let $(G,+)$ be a finite group.
 We denote the identity mapping on $G$ by $\1$, and
 let $\0: G\rightarrow G, x\mapsto 0$.
 A group $\Phi$ of automorphisms on $(G,+)$ is said to be
 \emph{fixed-point-free} (\emph{fpf}) iff $|\Phi(x)|=|\Phi|$ for all
 $x\in G^*$.

\begin{prop} \label{pro:sdffromfpf}
 Let $\Phi$ be group of fpf automorphisms on $(G,+)$, and let 
 $S\subseteq\Phi\cup\{\0\}$ with $|S|>1$. 
 Let $\mu,\nu$ be positive integers such that the following are satisfied:
%
\begin{enumerate}
\item \label{it:stab}
 $|G_{S(a)}| = \mu$ for all $a\in G^*$;
\item \label{it:sim}
 $|S(a)/\sim| = \nu$ for all $a\in G^*$.
\end{enumerate}
 Then $\{S(x)\ |\ x\in G^*\}$ is a
 $(|G|,|S|,\frac{|S|\cdot(|S|-1)}{\mu\nu})$-sdf.
\end{prop}

\begin{proof}
 Since $\Phi$ is fpf, 
 we have $|S(x)|=|S|$ for all $x\in G^*$.
 We note that for the same reason the map
 $\alpha-\beta: G\rightarrow G, x\mapsto \alpha(x)-\beta(x)$ is bijective
 for all $\alpha,\beta\in\Phi\cup\{\0\},\alpha\neq\beta$.
 Thus $(\alpha-\beta)^{-1}$ exists.

 It remains to check condition~\eqref{it:sdf4} of Definition~\ref{def:sdf}.
 Let $d\in G^*$, and let $\B=\{S(x)\ |\ x\in G^*\}$. Now we have
$$\begin{array}{ll}
    & \{(B,a,b)\ |\ B\in\B, a,b\in B, a-b=d \} = \\
  = & \{(S(x),\alpha(x),\beta(x))\ |\ x\in G^*, \alpha,\beta\in S, \alpha(x)-\beta(x)=d \} = \\
  = & \{(S(x),\alpha(x),\beta(x))\ |\ \alpha,\beta\in S, \alpha\neq\beta, x=(\alpha-\beta)^{-1}(d) \},
\end{array}$$
 and the cardinality of the last set is $|S|\cdot(|S|-1)$.
 Thus $\B$ is an sdf with parameters as given in the proposition.
\end{proof}

 
 Several of the designs from planar near-rings as described in 
 \cite{Clay:Nearrings} arise in the situation of
 Proposition~\ref{pro:sdffromfpf} for different choices of a set $S$ of
 endomorphisms on a group $G$.
 We will not prove the next 2 well-known results. We note that also for
 the original proofs the difficulty lies entirely in verifying the
 conditions~\eqref{it:stab} and~\eqref{it:sim} as given in
 Proposition~\ref{pro:sdffromfpf}.

\begin{cor} \cite[cf. p. 59, Theorem 5.5]{Clay:Nearrings}
 Let $\Phi$ be fpf on $G$. Then $\dev\{\Phi(x)\ |\ x\in G^*\}$ is a 
 $(|G|,|\Phi|,|\Phi|-1)$-design. 
\end{cor}

\begin{cor}
 Let $\Phi$ be fpf on $G$.
 Let $S=\Phi\cup\{0\}$, and let $\B=\{S(x)\ |\ x\in G^*\}$.
\begin{enumerate}
\item \cite[cf. p. 118, Theorem 7.9]{Clay:Nearrings}
 If $B$ is not a subgroup of $G$ for any $B\in\B$,
 then $\dev\B$ is a $(|G|,|\Phi|+1,|\Phi|+1)$-design.  
\item \cite[cf. p. 118, Theorem 7.11]{Clay:Nearrings}
 If $B$ is a subgroup of $G$ for all $B\in\B$,
 then $\dev\B$ is a $(|G|,|\Phi|+1,1)$-design. 
\end{enumerate}
\end{cor}

\begin{prop}  \label{transnormal}
 Let $\Phi$ be fpf on $G$, and let $S\subseteq\Phi\cup\{0\}$.
 Let $\Psi$ be a group of automorphisms of $G$, such that
 $\Psi$ normalizes $S$ and $\Psi$ is transitive on $G^*$.

 Then $\B=\{S(x)\ |\ x\in G^*\}$ is an sdf, and $G\rtimes\Psi$ acts as a
 doubly transitive group of automorphisms on $\dev\B$.
\end{prop}

\begin{proof}
 It suffices to check that the assumptions of Proposition~\ref{pro:sdffromfpf}
 are satisfied. 
 Let $a,b\in G^*, c\in G$ such that $S(a) = S(b) + c$, and let $\psi\in\Psi$.
 Since $\psi$ normalizes $S$, we have $S(\psi(a)) = S(\psi(b)) + \psi(c)$.
 Thus we obtain $|G_{S(\psi(a))}| = |G_{S(a)}|$ and
 $|S(\psi(a))/\sim| = |S(a)/\sim|$. Since $\Psi$ is transitive on $G^*$,
 the assumptions of Proposition~\ref{pro:sdffromfpf} are satisfied.
 The result follows.  
\end{proof}

\begin{cor}  \label{cor:nf}
 Let $(F,+,\cdot)$ be a finite left near-field, and let $S\subseteq F$
 such that $xS=Sx$ for all $x\in F$.
 Then $\B=\{Sx\ |\ x\in F^*\}$ is an sdf.
\end{cor}

\begin{proof} Straightforward by Proposition~\ref{transnormal}. \end{proof}
 
 We note that the assumption $xS=Sx$ for all $x\in F$ in the above corollary
 is trivially fulfilled if $(F,+,\cdot)$ is a field.

\section{A generalization of Sun's segments}

 We will now describe how the designs of segments as defined
 in~\cite{Sun:Segments} fit into our setting of short difference families.
 Our main result is the following.

\begin{prop} \label{pro:segments}
 Let $(G,+)$ be a group, and let $S$ be a set of endomorphisms on $G$ such that
 the following conditions are satisfied:
\begin{enumerate}
\item \label{it:Snontriv}
 $\0,\1\in S$ and $|S|>2$;
\item \label{it:S=1-S}
 $\{ \1-\alpha\ |\ \alpha\in S\} = S$; 
\item \label{it:Sfpf}
 $\langle S^*\rangle$ is a group of fpf automorphisms on $G$.
\end{enumerate}
 We assume that $|G|$ and $|\langle S^*\rangle|$ are odd.
 Then $\{S(x)\ |\ x\in G^*\}$ is a $(|G|,|S|,\frac{|S|\cdot(|S|-1)}{2})$-sdf.
\end{prop}

 Throughout this section we assume
 that $(G,+)$ is a group with a set of endomorphisms $S$ that satisfy the
 conditions~\eqref{it:Snontriv},~\eqref{it:S=1-S}, and~\eqref{it:Sfpf} of
 Proposition~\ref{pro:segments}.

 As an example of such a set $S$, we may consider a group of fpf automorphisms
 $\Phi$ on $G$ such that $6$ divides $|\Phi|$.
 Now let $S = \{\alpha\in\Phi\ |\ \ord\alpha = 6 \}\cup\{\0,\1\}$. 
 We show that if $\alpha\in\Phi$ has order $6$,
 then $\alpha^5 = \1-\alpha\in\Phi$. 

 Let $\alpha\in\Phi$ with $\ord\alpha = 6$. Then $\alpha^3$ is an fpf
 automorphism of order $2$. Thus $\alpha^3=-\1$ is an automorphism of $G$,
 and hence $G$ is abelian. Let $x\in G$.
 By $\alpha^2(x+\alpha^2(x)+\alpha^4(x))=\alpha^2(x)+\alpha^4(x)+x$
 we find that $x+\alpha^2(x)+\alpha^4(x)$ is a fixed-point of $\alpha^2\neq\1$.
 Thus $x+\alpha^2(x)+\alpha^4(x)=0$ for all $x\in G$.
 Now $\1+\alpha^2+\alpha^4=\0$ yields $\1+\alpha^4=-\alpha^2$.
 By $\alpha^3=-\1$, we obtain $\1-\alpha=\alpha^5$.
 

 For the proof of Proposition~\ref{pro:segments}, we need a bit of
 preparation.

\begin{lem} \label{le:Gisab} 
 Let $\alpha$ be an automorphism on $(G,+)$ such that $\1-\alpha$ is also
 an automorphism on $(G,+)$.
 Then $(G,+)$ is abelian.
\end{lem}

\begin{proof}
 Let $\beta = \1-\alpha$, and let $x,y\in G$. Then we have
\begin{eqnarray} \label{eq:hom1}
   \beta(x+y) = x+y-\alpha(x+y) = x+y-\alpha(y)-\alpha(x)
\end{eqnarray}
 and 
\begin{eqnarray} \label{eq:hom2}
  \beta(x+y) = \beta(x) + \beta(y) = x-\alpha(x)+y-\alpha(y). 
\end{eqnarray}
 By~\eqref{eq:hom1} and~\eqref{eq:hom2}, we obtain 
\begin{eqnarray} \label{eq:hom3} 
(y-\alpha(y))-\alpha(x) = -\alpha(x)+(y-\alpha(y)).
\end{eqnarray}
 Since $-\alpha$ and $\1-\alpha$ are bijections on $G$,
 equation~\eqref{eq:hom3} yields $a+b = b+a$ for all $a,b\in G$. 
\end{proof}

\begin{lem} \label{le:normaliscentral}
 Let $\Phi$ be fpf on $(G,+)$, and let $\varphi\in\Phi$.
 If $\{\alpha,\1-\alpha\}\subseteq N_\Phi(\langle\varphi\rangle)$, then
 $\{\alpha,\1-\alpha\}\subseteq C_\Phi(\langle\varphi\rangle)$.
\end{lem}

\begin{proof}
 Let $\{\alpha,\1-\alpha\}\subseteq N_\Phi(\langle\varphi\rangle)$.
 We note that $\langle\alpha,\1-\alpha\rangle$ is an abelian, and hence
 cyclic group of fpf automorphisms.
 Let $\beta\in\Phi$ generate $\langle\alpha,\1-\alpha\rangle$, and let
 $\alpha=\beta^i,\1-\alpha=\beta^j$.

 Since $\beta$ normalizes $\langle\varphi\rangle$, there exists an integer $r$
 such that $\varphi^\beta = \varphi^r$.
 Then $\varphi^\alpha = \varphi^{(r^i)}$, and
 $\varphi^{\1-\alpha} = \varphi^{(r^j)}$.

 By
 $\varphi(\1-\alpha) = \varphi-\varphi\alpha = \varphi-\alpha\varphi^{(r^i)}$,
 we have $(\1-\alpha)\varphi^{(r^j)} = \varphi-\alpha\varphi^{(r^i)}$.
 Thus we obtain
\begin{eqnarray} \label{eq:varphi}
 \alpha(\varphi^{(r^i)}-\varphi^{(r^j)}) & = & \varphi-\varphi^{(r^j)}.
\end{eqnarray}

 If $\varphi^{(r^i)}\neq\varphi^{(r^j)}$, then
 $\varphi^{(r^i)}-\varphi^{(r^j)}$ is invertible because $\Phi$ is fpf.
 Then
 $\alpha=(\varphi-\varphi^{(r^j)})(\varphi^{(r^i)}-\varphi^{(r^j)})^{-1}$,
 and, in particular, $\alpha$ commutes with $\varphi$.
 But then also $\1-\alpha$ commutes with $\varphi$. We obtain
 $\varphi^{(r^i)}=\varphi=\varphi^{(r^j)}$ in contradiction to our assumption.
 Thus we have $\varphi^{(r^i)}=\varphi^{(r^j)}$, and by (\ref{eq:varphi}) 
 $\varphi=\varphi^{(r^j)}$. Both $\alpha$ and $\1-\alpha$ commute with
 $\varphi$. 
\end{proof}

\begin{cor} \label{cor:mciscyc}
 If $\langle S^*\rangle$ is metacyclic, then $\langle S^*\rangle$ is cyclic.
\end{cor}

\begin{proof}
 Let $\Phi = \langle S^*\rangle$ be a metacyclic group of fpf automorphisms.
 Then there exists $\varphi\in\Phi$ such that $N = \langle\varphi\rangle$
 is normal in $\Phi$, and there exists $\psi\in\Phi$ such that
 $\psi N$ generates $\Phi/N$.
 
 By Lemma~\ref{le:normaliscentral}, a cyclic normal subgroup of $\Phi$ is
 central in $\Phi$. In particular, $\psi$ and $\varphi$ commute.
 Thus $\Phi$ is abelian and hence cyclic.
\end{proof}

 We note that, by Corollary~\ref{cor:mciscyc}, $\Phi=\langle S^*\rangle$ is 
 cyclic if $|\Phi|$ is odd.

\begin{cor} \label{cor:classification}
 Let $S$ be a set of endomorphisms of $(G,+)$ such that $S=\1-S$ and
 $\Phi = \langle S^*\rangle$ is a group of fpf automorphisms of $G$.

 Then $\Phi/Z(\Phi)$ is trivial or isomorphic to one of the following:
 $A_4$, $S_4$, $A_5$, $S_5$.
\end{cor}

\begin{proof}
 By~\cite[Theorem 1.4]{Brown:Frobenius}, we have  a unique normal subgroup $N$
 of $\Phi$ such that all Sylow subgroups of $N$ are cyclic and $\Phi/N$
 is isomorphic to one of the following groups: $1, V_4,  A_4, S_4, A_5, S_5$ 
 where $1$ denotes the trivial group and $V_4$ denotes the Klein group.
 We note that $N'$ is cyclic and a normal subgroup of $\Phi$.
 Thus the elements of $S^*$ normalize $N'$ and,
 by Lemma~\ref{le:normaliscentral}, $S^*$ centralizes $N'$.
 Hence we have $N' \subseteq Z(\Phi)$.
 In particular, $N$ is nilpotent. Now $N$ is the direct product of its cyclic
 Sylow subgroups, and $N$ is cyclic. As a cyclic normal subgroup, $N$ is 
 central in $\Phi$ by Lemma~\ref{le:normaliscentral}.
 Thus $N \subseteq Z(\Phi)$.

 First we assume that $\Phi/N$ is abelian. Since $N$ is central, this yields
 that $\Phi$ is nilpotent. Let $P$ denote the Sylow $2$-subgroup of $\Phi$. 
 Then $P$ is a cyclic group or a generalized quaternion group. In any case,
 $P$ has a normal cyclic subgroup $R$ of index $2$. Since $\Phi$ is the
 direct product of its Sylow subgroups, $R$ is normal in $\Phi$.
 By Lemma~\ref{le:normaliscentral}, $R$ is central in $\Phi$, which implies
 that $P$ is cyclic. Thus $\Phi/N$ is not isomorphic to $V_4$.

 It remains that $\Phi/N$ is isomorphic to $1, A_4, S_4, A_5$, or $S_5$.
 Since all these groups have trivial center, we finally obtain $N = Z(\Phi)$.
 The corollary is proved.
\end{proof}

\begin{lem}  \label{le:eqstab}
 Let $a,b\in G^*,c\in G$. If $S(a)=S(b)+c$, then $-a+b+2c \in G_{S(a)}$.
\end{lem}

\begin{proof}
 Let $S(a)=S(b)+c$. Since $(G,+)$ is abelian and $S(x) = x-S(x)$ for all
 $x\in G^*$, we obtain:
\begin{eqnarray*}
 S(a)-a+b+2c & = & -S(a)+b+2c \\
           & = & -(S(b)+c)+b+2c \\
           & = & -S(b)+b+c \\
           & = & S(b)+c \\
           & = & S(a)
\end{eqnarray*}
 Thus $-a+b+2c \in G_{S(a)}$.
\end{proof}


\begin{lem}  \label{le:eqtrans}
 Let $\Phi=\langle S\setminus\{\0\}\rangle$. We assume that $|\Phi|$ is odd
 and that $G_{S(a)}=\{0\}$ for all $a\in G^*$.
 Then $S(a)=S(b)+c$ iff ($b=a$ and $c=0$) or ($b=-a$ and $c=a$).
\end{lem}

\begin{proof}
 Let $S(a)=S(b)+c$. By Lemma~\ref{le:eqstab}, we have $-a+b+2c \in G_{S(a)}$.
 Since $G_{S(a)} = \{0\}$ by assumption, we obtain $b=a-2c$ and
 $S(a)=S(a-2c)+c$.

 We note that $S(a-2c)$ contains the elements $0$ and $a-2c$.
 By $S(a)=S(a-2c)+c$, we obtain that $0+c$ and $(a-2c)+c$ are in $S(a)$.
 With $\{0,a\}\subseteq S(a)$ and $S(a)-c=S(a-2c)$ we find
 $\{-c,a-c\}\subseteq S(a-2c)$. 
 Summing up, we have $\{0,a,c,a-c\}\subseteq S(a)$ and
 $\{0,a-2c,-c,a-c\}\subseteq S(a-2c)$.

 Seeking a contradiction, we suppose that $c\neq a$ and $c\neq 0$.
 Let $S^* = S\setminus\{\0\}$.
 Since $a-c\in S^*(a)$ and $a-c\in S^*(a-2c)$, both $S^*(a)$ and $S^*(a-2c)$ 
 are contained in $\Phi(a-c)$. In particular, $c$ and $-c$ are in $\Phi(a-c)$.

 If $c=-c$, then $2c=0$ and $S(a)=S(a)+c$. Now $G_{S(a)}=0$ yields $c=0$,
 which contradicts our assumption that $c\neq 0$.
 
 If $c\neq-c$, then there exists $\varphi\in\Phi,\varphi\neq\1$, such that
 $\varphi(c)=-c$. Now, $c\neq 0$ is a fixed-point of $\varphi^2$.
 Hence $\varphi^2=\1$ and $|\Phi|$ is even, which contradicts the assumption
 of the lemma. Thus we have $c=a$ or $c=0$. The lemma is proved.
\end{proof}

\emph{Proof of Proposition~\ref{pro:segments}.}
 Let $G$ and $S$ satisfy the assumptions, and let $a\in G^*$.
 We will prove that 
\begin{equation} \label{eq:stabistriv}
 |G_{S(a)}| = 1. 
\end{equation}
 Seeking a contradiction, we suppose that there exists $d\in G_{S(a)}$ such
 that $d\neq 0$. Since $\0\in S$, we then have $\{d,-d\}\in S(a)$.
 For $\Phi = \langle S^*\rangle$, we obtain $\{d,-d\}\in\Phi(a)$.
 In particular, there is $\varphi\in\Phi$ such that $\varphi(d) = -d$.
 We have $d\neq -d$ and $\varphi\neq\1$ by the assumption that $|G|$ is odd.
 Since $\varphi^2$ fixes $d\neq 0$, we then have $\ord\varphi = 2$.
 This contradicts the assumption that $|\Phi|$ is odd. 
 Thus we have $d = 0$, and~\eqref{eq:stabistriv} is proved.
 Next we show that
\begin{equation} \label{eq:simistriv}
 |S(a)/\sim| = 2. 
\end{equation}
 For $b\in G^*$,$c\in G$ such that $S(a)=S(b)+c$, Lemma~\ref{le:eqtrans}
 yields that $b=a$ or $b=-a$. Thus we obtain $S(a)/\sim = \{S(a),S(-a)\}$.
 Since $a\neq -a$ by the assumption that $|G|$ is odd, we have
 $S(a)\neq S(-a)$. This proves~\eqref{eq:simistriv}.

 Now we obtain the result from~\eqref{eq:stabistriv} and~\eqref{eq:simistriv}
 with Proposition~\ref{pro:sdffromfpf}.  
\fixqed


 We give several examples for designs of segments from short difference
 families.
 
\begin{exa}
 If $G$ is an elementary abelian $2$-group, then we have $S = 1-S = 1+S$.
 Thus we know that $\{0,a\}\subseteq G_{S(a)}$ for all $a \in G^*$.
 Assume equality, that is, $G_{S(a)} = \{0,a\}$ for all $a \in G^*$.
 Then Lemma~\ref{le:eqstab} yields $|S(a)/\sim| = 1$ for all $a \in G^*$.
 By Proposition~\ref{pro:sdffromfpf}, $\{S(a)\ |\ a\in G^*\}$ forms a
 $(|G|,|S|,\frac{|S|\cdot(|S|-1)}{2})$-sdf for $G$.
\end{exa}

\begin{exa}
 Let $|S| = 3$, that is, $S=\{\0,\1,\2^{-1}\}$ where $\2^{-1}$ is the inverse
 mapping of $\2: G\rightarrow G, x\mapsto x+x$.
 Suppose there is $a\in G^*$ such that $G_{S(a)} \neq \{0\}$. 
 Then $G_{S(a)} = S(a)$, which yields $\2^{-1}(a) = -a$.
 Hence $\2^{-1}$ is an fpf automorphism of order $2$, and $\2^{-1}= \2 = -\1$.
 Now $2x=-x$ for all $x\in G^*$ yields $3x=0$; $G$ is an elementary
 abelian $3$-group.
 Thus $S(x) = \langle x\rangle = G_{S(x)}$ for all $x\in G$. 
 By Proposition~\ref{pro:sdffromfpf}, $\{S(a)\ |\ a\in G^*\}$ forms a
 $(|G|,3,1)$-sdf for $G$.
 
 If $|S| = 3$ and $G_S(a) = \{0\}$ for all $a\in G^*$, then
 $S(a)/\sim = \{S(a),-S(a)\}$ for all $a\in G^*$ by Lemma~\ref{le:eqstab},
 and $\{S(a)\ |\ a\in G^*\}$ forms a $(|G|,3,3)$-sdf for $G$ by
 Proposition~\ref{pro:sdffromfpf}.
\end{exa}

\begin{exa}
 Let $S \subseteq \{\alpha\in\Phi\ |\ \ord\alpha = 6 \}\cup\{\0,\1\}$ such
 that $\{\0,\1\} \in S$ and $\alpha^5 \in S$ for all $\alpha \in S$.
 Then $|S|$ is even, $|S| \leq 22$, and $\{S(a)\ |\ a\in G^*\}$ is a
 $(|G|,|S|,\frac{|S|\cdot(|S|-1)}{2})$-sdf for $G$.

 
 From Corollary~\ref{cor:classification} we obtain that $\langle S^* \rangle$
 is cyclic, isomorphic to $\SL(2,3)$ or to $\SL(2,5)$. Thus $G$ has at most
 $20$ fixed-point-free automorphisms of order $6$, and $|S| \leq 22$.
 
 Next we show that $G_{S(a)} = \{0\}$ for all $a\in G^*$. Let $a\in G^*$,
 $c\in G_{S(a)}$. Then we have $\alpha\in S$ such that $c = \alpha(a)$.
 Since $-c \in G_{S(a)}$, we also have $\beta\in S$ such that $-c = \beta(a)$.
 We note that $|G|$ is odd because $\langle S^* \rangle$ is even.
 Suppose that $\alpha\neq\0$. Then $c\neq 0$ and
 $\beta\alpha^{-1}(c) = -c \neq c$.
 By assumption $\beta\alpha^{-1}$ is fixed-point-free and has order $2$.
 Hence $\beta\alpha^{-1} = -\1$ and $\beta = -\alpha$. If $\alpha = \1$,
 then we obtain $\beta = -\1\in S$ which contradicts the assumption on $S$.
 If $\ord\alpha = 6$, then $\beta = \alpha^4 \in S$ and $\ord\beta = 3$
 which yields a contradiction. Thus we have $\alpha = 0$ and
 $G_{S(a)} = \{0\}$.

 Next we show that $S(a)/\sim = \{S(a),-S(a)\}$ for all $a\in G^*$.
 Let $a,b\in G^*, c\in G$ such that $S(a) = S(b) + c$. Seeking a contradiction
 we suppose that $c\not\in\{0,a\}$.
 Then we have $\alpha\in S\setminus\{\0,\1\}$ and $\beta\in S^*$ such that
 $c = \alpha(a)$ and $c = -\beta(b)$. We note that $\beta\neq\1$ because
 otherwise $S(a) = S(-c) +c = S(c)$ and Lemma~\ref{le:eqstab} yields $a=c$.
 Hence both $\alpha$ and $\beta$ have order $6$. From Lemma~\ref{le:eqstab}
 we obtain $-\alpha^{-1}(c) - \beta^{-1}(c) + 2c = 0$.
 Together with $\alpha^{-1} = \1-\alpha$ and $\beta^{-1} = \1-\beta$
 (see the remark below Proposition~\ref{pro:segments}),
 this yields $\alpha(c) + \beta(c) = 0$. Thus $\beta = \alpha^4$ by the 
 assumption that $\langle S^* \rangle$ is fpf. This contradicts that 
 $\alpha,\beta$ have order $6$. Hence $\alpha,\beta \in \{\0,\1\}$
 We have either $c = 0$, $b = a$ or $c = a, b = -a$. 
 Thus $S(a)/\sim = \{S(a),-S(a)\}$. Proposition~\ref{pro:sdffromfpf} yields
 that $\{S(a)\ |\ a\in G^*\}$ is an sdf.
\end{exa}

\begin{exa} 
 We note that for a left near-field $(F,+,\cdot)$ each element $a\in F^*$
 induces an automorphism $\lambda_a: x\mapsto a\cdot x$ on $(F,+)$.
 Let $\Phi=\{\lambda_a\ |\ a\in F^*\}$. Then $\Phi$ is a group of fpf
 automorphisms of  $(F,+)$, and $\Phi$ is isomorphic to $(F^*,\cdot)$. 

 With $1$ denoting the identity of $(F,+,\cdot)$ we have $\lambda_1 = \1$.
 If $(F,+,\cdot)$ is not a field, then we do not have
 $\1-\lambda_a = \lambda_{1-a}$ in general.
 While $\1-\lambda_a$ for $a\neq 1$ is an automorphism of $(F,+)$, in general
 it is not true that $\1-\lambda_a\in\Phi$. We note that
 \[ \Phi\cap(\1-\Phi) = \{\lambda_a\ |\ a\in F^* \text{ such that }(1-a)x=x-ax\text{ for all }x\in F\}.\]  

 Let $(F,+,\cdot)$ be a Dickson near-field with Dickson pair $(q,n)$.
 That is $(F,+,\cdot)$ is coupled to the field $\GF(q^n)=(F,+,*)$.
 Let $\Phi = \{ \lambda_f\ |\ f\in F^*\}$. Then there
 exist $a,b\in F$ such that $\Phi$ is generated by
 $\alpha: x\mapsto a*x$ and $\beta: x\mapsto b*x^{q}$.
 That is, $\Phi$, and hence $(F^*, \cdot)$ is metacyclic.
 Let $S^*\subseteq\Phi$ and $S=\1-S$. By Lemma~\ref{le:normaliscentral}
 we have that $\langle S^*\rangle \subseteq \langle\alpha\rangle$, and 
 in particular, $\langle S^*\rangle$ is cyclic.
 All mappings in $S$ are of the form $\mu_f: x\mapsto f*x$ where
 $f\in\GF(q^n)$ and $*$ is the multiplication in the field $\GF(q^n)$.

 Let $T=\{ f\in F\ |\ \mu_f\in S\}$. Then $S(x) = T*x$ for all $x\in F$. 
 Since $T$ is central in the multiplicative group of $\GF(q^n)$,
 Corollary~\ref{cor:nf} yields that $\{T*x\ |\ x\in F^*\}$ is an sdf. 
\end{exa}

\section{Conclusion and further work}

This paper has presented an explanation of the creation of
designs from planar near-rings as pioneered by Ferrero and
Clay, using a generalization of the well known concept of difference families.
A further analysis of work with segments was possible,
both in its geometric view (the definition of segments by their
endpoints) and the design view (determination of
some of the resulting designs).

No clear understanding of the relation of circularity \cite{beidarFongKe96,keThesis,kekiechle96}
to this structure has been seen. Future work may investigate this, as
it seems to relate strongly to the motivating examples of 
planar near-rings from the complex plane, as does the
segments construction.

\bibliographystyle{alpha}
\bibliography{project}

\end{document}